\documentclass[12pt]{article} 
\usepackage{amssymb,amsfonts,mathrsfs,graphicx}

\usepackage[utf8]{inputenc}   

\usepackage{abstract}

\usepackage{color}

\definecolor{dmagenta}{rgb}{0.4,0,0.4}   
\definecolor{dblue}{rgb}{0,0,0.4}
\definecolor{dgreen}{rgb}{0,0.5,0}
\definecolor{dred}{rgb}{0.6,0,0}

\newtheorem{theorem}{Theorem}[section]

\newtheorem{proposition}[theorem]{Proposition}

\newtheorem{definition}[theorem]{Definition}

\newtheorem{rmrk}[theorem]{Remark}

\makeatletter
\@addtoreset{equation}{section}
\makeatother

\newcommand{\fig}[3] {
\medskip\smallskip
\begin{figure}[htb]
  \centering
  \includegraphics[width=#2]{#1.pdf}
  \begin{minipage}[t]{0.80\linewidth} 
    \caption{#3}
    \protect\label{#1}
  \end{minipage}
\end{figure}
\medskip
}

\newcommand{\R} {\mathbb{R}}

\newcommand{\Z} {\mathbb{Z}}
\newcommand{\N} {\mathbb{N}}
\newcommand{\qed} {\hfill {\small Q.E.D.} \par\medskip}
\newcommand{\skippar} {\par\medskip}

\newcommand{\proof} {\noindent \textsc{Proof.} }

\newcommand{\article}[3] {\textsc{{#1}}, {\itshape {#2}}, {{#3}}.}
\newcommand{\book}[3] {\textsc{{#1}}, {\itshape {#2}}, {{#3}}.}
\newcommand{\vol} {\textbf}

\newcommand{\rset}[2] {\left\{ #1 \: \left| \: #2 \right. \! \right\} }

\newcommand{\symmdiff} {\triangle}

\newcommand{\om} {\omega}
\newcommand{\gm} {\gamma}       
\renewcommand{\l} {\ell}

\newcommand{\mN} {\mathcal{N}}
\newcommand{\tu} {\mathcal{T}}  
\newcommand{\ta} {\mathcal{Q}}  
\newcommand{\mO} {\mathcal{O}}

\newcommand{\mM} {\mathcal{M}}

\newcommand{\bi} {billiard}
\newcommand{\fn} {function}
\newcommand{\me} {measure}
\newcommand{\tr} {trajector}
\newcommand{\erg} {ergodic}
\newcommand{\sy} {system}
\newcommand{\hyp} {hyperbolic}
\newcommand{\pr} {probability}

\newcommand{\dsy} {dynamical system}
\renewcommand{\o} {orbit}

\begin{document}

\title{\textbf{Recurrence for quenched random Lorentz tubes}}

\author{\textsc{Giampaolo Cristadoro$^1$, Marco Lenci$^1$, 
Marcello Seri$^{1,2}$}
\vspace{15pt}
\\ 
$^1$ Dipartimento di Matematica, Universit\`a di Bologna, 
\\ 
Piazza di Porta San Donato 5, 40126 Bologna, Italy 
\vspace{10pt}
\\ 
$^2$ Department of Mathematics, University of Erlangen-Nuremberg, 
\\ 
Bismarckstr.\ 1 1/2, 91053 Erlangen, Germany
\vspace{5pt}
}

\date{Version published on Chaos \vol{20} (2010), 023115 \\ 
+ correction of erratum \\ 
\vspace{10pt} November 2010}

\maketitle

\begin{abstract}
  We consider the billiard dynamics in a strip-like set that is
  tessellated by countably many translated copies of the same
  polygon. A random configuration of semidispersing scatterers is
  placed in each copy. The ensemble of dynamical systems thus defined,
  one for each global choice of scatterers, is called \emph{quenched
  random Lorentz tube}. We prove that, under general conditions,
  almost every system in the ensemble is recurrent.

\end{abstract}

\bigskip

{\bfseries 
  A Lorentz tube is a system of a particle (or, from a statistical
  viewpoint, many non-interacting particles) freely moving in a domain
  extended in one direction and performing elastic collisions with
  randomly placed obstacles. These kinds of ``extended billiards''
  are, on the one hand, paradigms of systems where some transport
  properties can be studied in a rigorous mathematical way and, on the
  other hand, reliable models for real situations, such as transport
  in nanotubes, heat diffusion and molecular dynamics in wires or
  other disordered tubular settings, etc. The primary interest in
  their study lies on such properties as recurrence, diffusivity, and
  transmission rates. Unfortunately, few rigorous results are
  available and their proofs typically rely on some periodic
  structure. In this paper a more realistic situation is taken into
  account: the so-called quenched disorder. Recurrence is proved for
  almost every realization of the configuration of obstacles,
  impliying strong chaotic properties for these types of systems.  
}

\section{Introduction}
\label{sec-int}

This paper concerns the dynamics of a particle in certain
two-dimensional \sy s which are infinitely extended in one
dimension. More precisely, we will study \dsy s in which a point
particle moves in a strip (or similar set) $\tu \subset \R^2$, which
contains a countable number of convex scatterers, see the example in
Fig.~\ref{imgquadrati1}.  The motion of the particle is free until it
collides with either the boundary of $\tu$ or a scatterer, both of
which are thought to have infinite mass. The collisions are totally
elastic, so they obey the usual Fresnel law: the angle of reflection
equals the angle of incidence.

\fig{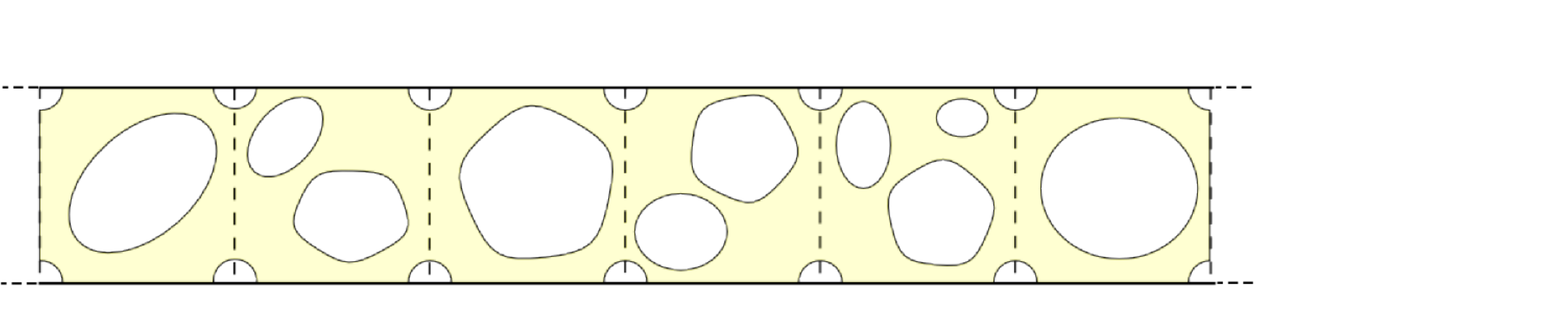}{0.9\linewidth}{A simple Lorentz tube.}

In the taxonomy of \dsy s, these models belong to the class of
semidispersing planar \bi s. In particular, they are extended
semidispersing \bi s, which very much resemble a Lorentz gas. We thus
call them \emph{effectively one-dimensional Lorentz gases} or, more
concisely, \emph{Lorentz tubes (LTs)}.

Systems like these (especially their three-dimensional counterparts,
cf.\ last paragraph of this section) find application in the sciences as
models for the dynamics of particles (e.g., gas molecules) in narrow
tubes (e.g., carbon nanotubes). A very minimal list of references,
from the more experimental to the more mathematical, includes
\cite{heal}, \cite{acm}, \cite{lwwz}, \cite{cmp}, \cite{aacg},
\cite{fy}, \cite{f}. (See further references in those papers.) An
interesting fact is that both experimentalists and theoreticians seem
to have a primary interest --- sometimes for different reasons --- in
the same question, namely the diffusion properties of these gases. As
we discuss below, this is our case as well, although the results we
present in this note must be considered preliminary in this respect.

From a mathematical viewpoint, LTs are interesting because they are
among the very few extended \dsy s, with a certain degree of realism,
that mathematicians can prove something about. By the ill-defined
expression \emph{extended \dsy} we generally mean a \dsy\ on a
non-compact phase space whose physically relevant (invariant) \me\ is
infinite. For such \sy s, the very fundamentals of ordinary \erg\
theory do not work \cite{a}: for example, the Poincar\'e Recurrence
Theorem fails to hold and one does not know whether the \sy\ is totally
recurrent (almost every point returns arbitrarily close to its initial
condition), totally transient (almost every point escapes to
infinity), or mixed.

In fact, as it turns out, recurrence is not just the most basic
property one wants to establish in order to even consider studying the
chaotic features of an extended \dsy\ (it is sometimes said that, if
\erg ity is the first of a whole hierarchy of stochastic properties
that a \dsy\ can possess, recurrence is the \emph{zeroeth property}); 
for a Lorentz gas at least, a number of stronger \erg\ properties follow
from recurrence: for example, \erg ity of the extended \dsy, $K$-mixing
of the first-return map to a given scatterer, etc.\ \cite{l1}.

Let us briefly explain our model. We consider the connected set $\tu
\subset \R^2$ tessellated by the repetition, under the action of $\Z$,
of a given \emph{fundamental domain} $C$, which we assume to be a
polygon.  In each copy of $C$, henceforth referred to as \emph{cell},
we place a random configuration of convex \emph{scatterers}, according
to some rule that we specify later. Given a global configuration of
scatterers, we consider the \bi\ dynamics in the complement (to $\tu$)
of the union of all the scatterers.

So, each model just described does not correspond to one \dsy, but to
an \emph{ensemble} of \dsy s. In other words, we have a \emph{quenched
random \dsy}, in the sense that first a \sy\ is picked from a random
family and then its (deterministic) dynamics is observed. This
contrasts with random \dsy s, such as the random \bi\ channels of
\cite{fy}, \cite{f}, in which a new random map is applied at every
iteration of the dynamics.

Quenched random LTs are a bit more realistic and understandably harder
to study than random LTs, which are in turn harder than
\emph{periodic} LTs (when the configuration of scatterers is the same
in every cell). The same can be said of Lorentz gases which are
infinitely extended in both dimensions \cite{l2}. In fact, while
recurrence, the Central Limit Theorem (CLT) and several strong
stochastic properties are known for periodic Lorentz gases --- at
least under the so-called \emph{finite horizon} condition --- very
little is known for random or quenched random Lorentz gases (although
results were established for toy versions: \cite{l3}, \cite{als},
\cite{l4}).

As it turns out, when the effective dimension $\nu$ equals 2,
recurrence and the CLT go hand in hand, as a remarkable theorem by
Schmidt (Theorem \ref{thm-co-rec} below) shows \cite{s, l2}. This
provides another strong motivation for the study of the diffusive
properties of these gases, cf.\ also \cite{cd}.

We state the paper's main result in plain English, leaving a more
rigorous description to the remainder of the article, in particular
Section \ref{sec-povp}.  

This paper's main result is the \emph{almost sure} recurrence of our
quenched random LTs, under very mild geometrical conditions which
include the finite-horizon condition. Almost sure recurrence means
that almost every LT in the ensemble is Poincar\'e recurrent. To our
knowledge, this is the first time that recurrence is proved for the
typical element of a fairly general class of Lorentz gases (albeit
effectively one-dimensional). The main ingredient of the proof is the
above-mentioned theorem by Schmidt, which is particularly powerful for
$\nu = 1$.

\medskip

The exposition is organized as follows: In Section \ref{sec-pre} we
give a precise definition of our LTs and state some of their
properties. Then in Section \ref{sec-rec} we introduce the tools that
we use to prove almost sure recurrence, namely Schmidt's Theorem and
an \erg\ \dsy\ endowed with a suitable one-dimensonal cocycle. The
latter objects are presented in detail in Section \ref{sec-povp},
where the main proof of the article is also given. Finally, in
Sections \ref{sec-ext}, we discuss some generalizations of our result.

Due to its technicality and lesser strength, the very important
generalization to the higher dimensional case will be presented
elsewhere.

\bigskip
\noindent
\textbf{Acknowledgments.}\ We thank Gianluigi Del Magno and Nikolai
Chernov for some illuminating discussions.

\section{Preliminaries and main assumptions}
\label{sec-pre}

We present the \sy\ in detail. Let $C_0$ be a closed polygon embedded
in $\R^2$, such that two of its sides, denoted $G^1$ and $G^2$, are
parallel and congruent. Then call $\tau$ the translation of $\R^2$
that takes $G^1$ into $G^2$, and define $C_n := \tau^n (C_0)$, with $n
\in \Z$. Each $C_n$ is called a \emph{cell} and $\tu := \bigcup_{n \in
\Z} C_n$ is called the \emph{tube}, see
Figs.~\ref{imgquadrati1}-\ref{imgtubostrano1}.

\fig{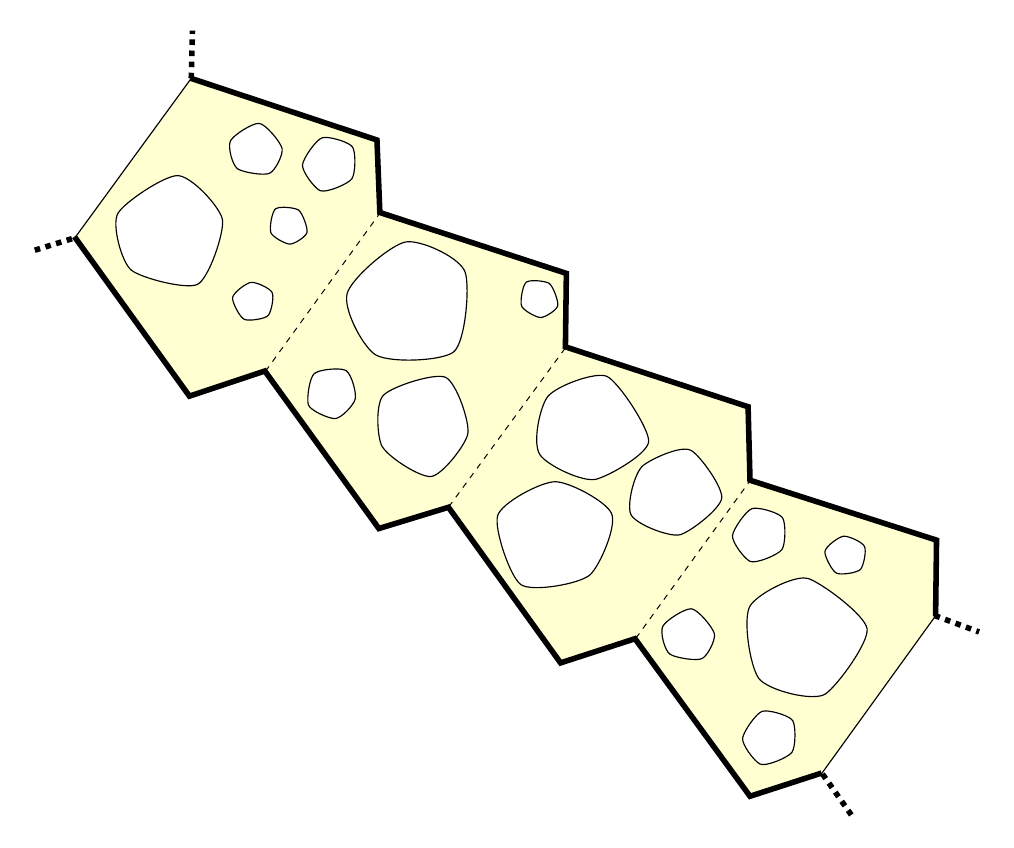}{0.6\linewidth}{A less trivial Lorentz tube.}

In every cell $C_n$ there is a configuration of closed, pairwise
disjoint, piecewise smooth, convex sets $\mO_{n,i} \subset C_n$ ($i =
1, \ldots, N$) which we call \emph{scatterers}. (Note that some
$\mO_{n,i}$ might be empty, so different cells might have a different
number of scatterers.) Each $\mO_{n,i} = \mO_{n,i}(\ell_n)$ is indeed
a \fn\ of the random parameter $\ell_n \in \Omega$, where $\Omega$ is
a \me\ space whose nature is irrelevant. The sequence $\ell := (
\ell_n )_{n \in \Z} \in \Omega^\Z$, which thus describes the global
configuration of scatterers in the tube $\tu$, is a stochastic process
obeying the \pr\ law $\Pi$. We assume that
\begin{itemize}
\item[(A1)] $\Pi$ is \erg\ for the left shift $\sigma: \Omega^\Z
  \longrightarrow \Omega^\Z$.
\end{itemize}

For each realization $\ell$ of the process, we consider the \bi\ in
the \emph{table} $\ta_\ell := \tu \setminus \bigcup_{n \in \Z}
\bigcup_{i=1}^N \mO_{n,i}(\ell_n)$. This is the \dsy\ $(\ta_\ell
\times S^{1}, \phi_\ell^t, m_\ell)$, where $S^{1}$ is the unit circle
in $\R^2$ and $\phi_\ell^t : \ta_\ell \times S^{1} \longrightarrow
\ta_\ell \times S^{1}$ is the \emph{\bi\ flow}, whereby $(q_t, v_t) =
\phi_\ell^t (q, v)$ represents the position and velocity at time $t$ of
a point particle with initial conditions $(q, v)$, undergoing free
motion in the interior of $\ta_\ell$ and Fresnel collisions at
$\partial \ta_\ell$.  (Notice that in this Hamiltonian \sy\ the
conservation of energy corresponds to the conservation of speed, which
is thus conventionally fixed to 1.)

Evidently, the above definition is a bit ambiguous since $\phi_\ell^t$
is discontinuous and there is a set of initial conditions for which it
is not even well defined. We thus declare that $t \mapsto \phi_\ell^t$
is right-continuous (i.e., if $t$ is a collision time, $v_t$ is the
\emph{post}-collisional velocity) and that a material point that hits
a non-smooth part of $\partial \ta_\ell$ stays trapped there forever
(assumption (A2) below ensures that this can only happen to a
negligible set of \tr ies).

Finally, $m_\ell$ is the Liouville invariant \me\ which, as is well
known, is the product of the Lebesgue \me\ on $\ta_\ell$ and the Haar
\me\ on $S^{1}$.

We call this \sy\ the \emph{LT corresponding to the realization
$\ell$}, or simply the \emph{LT $\ell$}. In the reminder, whenever
there is no risk of ambiguity, we drop the dependence on $\ell$ from
all the notation.

\skippar

The following are our assumptions on the geometry of the LT:

\begin{itemize}
\item[(A2)] There exist a positive integer $K$ such that, for
  $\Pi$-a.e.\ realization $\ell \in \Omega^\Z$, $\partial \mO_{n,i}$
  is made up of at most $K$ compact connected $C^3$ pieces, which may
  intersect only at their endpoints. These points will be referred to
  as \emph{vertices}.
\end{itemize}
Denoting, as we will do throughout the paper, $x := (q,v)$, let
$\gm(x)$ be the first time at which the point with initial conditions
$x$ hits a \emph{non-flat} part of the boundary (so this is not
exactly the usual free flight \fn!). Also, if $q$ is a smooth point of
$\partial \ta$, let $k(q)$ be the curvature of $\partial \ta$ at
$q$. We have:
\begin{itemize}
\item[(A3)] There exist two positive constants $\gm_m < \gm_M$ such
  that, for a.e.\ $\ell$ and all $x = (q,v)$ with $q \in \partial
  \ta$ and $k(q)>0$,
  \begin{displaymath}
    \gm_m \le \gm(x) \le \gm_M.
  \end{displaymath}
  Also, starting from any such $x$ and within the time $\gm(x)$, there
  cannot be more than $M$ collisions with flat parts of the boundary,
  where $M$ is a universal constant.

\item[(A4)] There exists a positive constant $k_m$ such that, for
  a.e.\ $\ell$, given a smooth point $q$ of the boundary, either
  $\partial \ta$ is totally flat at $q$ or
  \begin{displaymath}
    k(q) \ge k_m. 
  \end{displaymath}
\end{itemize}
In the language of \bi s, a \emph{singular \tr y} is a \tr y which, at
some time, hits the boundary of the table tangentially or in a
vertex. It follows that a finite segment of a \emph{non}-singular \tr
y depends continuously on its initial condition.  Also notice that, by
(A2), the set of all singular \tr ies is a countable union of smooth
curves in $\ta \times S^1$ and thus has \me\ zero. The next assumption
is meant to exclude pathological situations:
\begin{itemize}
\item[(A5)] For a.e.\ $\ell$ and all $i,j \in \{ 1,2 \}$, there is a
  non-singular \tr y entering $C_0$ through $G^i$ and leaving it
  through $G^j$.
\end{itemize}

A convenient way to represent a continuous-time \dsy\ is to select a
suitable Poincar\'e section and consider the first-return map
there. For \bi s, the section is customarily taken to be the set of
all pairs $(q,v) \in \partial \ta \times S^{1}$, where $v$ is a
post-collisional unit vector at $q$ (hence an inner vector relative to
$\ta$). Here we slightly modify this choice.

For $n \in \Z$ and $j \in \{ 1,2 \}$, denote by $G_n^j := \tau^n
(G^j)$ the side of $C_n$ corresponding to $G^j$ in $C_0$ ($G_n^1$ and
$G_n^2$ may be called the \emph{gates} of $C_n$, whence the notation).
Let $o_j$ be the inner normal to $G_n^j$, relative to $C_n$. Notice
that, under our hypotheses, $o_2 = -o_1$. Define
\begin{equation}
  \label{def-nnj}
  \mN_n^j : = \rset{(q,v) \in G_n^j \times S^1} {v \cdot
    o_j > 0}.
\end{equation}
The cross section we use is
\begin{equation}
  \label{def-nbar}
  \mM := \bigcup_{n \in \Z} \, \bigcup_{j=1,2} \mN_n^j,
\end{equation}
whose corresponding Poincar\'e map we denote $T = T_\ell$. In other
words, we only consider those times at which the particle crosses one
of the gates. In the lingo of \bi s, cross sections like these are
sometimes called ``transparent walls''. The Liouville \me\ for the
flow induces on a transparent wall an invariant \me\ given by
$d\mu(q,v) = (v \cdot o_q) \, dq dv$, where $o_q$ is the normal to the
section at $q$, directed towards the outgoing side of $(q,v)$
\cite{cm} (in our case, $o_q = o_j$ whenever $q \in \mN_n^j$).

So we end up with the \dsy\ $(\mM, T_\ell, \mu)$, whose invariant \me\
is infinite and $\sigma$-finite. Notice that, by design, the only
object that depends on the random configuration is the map $T_\ell$.

\skippar

In order to discuss the \hyp\ properties of this \sy, we need to
introduce its local stable and unstable manifolds (LSUMs). Since our
exposition does not require a rigorous definition of these objects, we
shall refrain from providing one, and point the interested reader to
the existing literature, e.g., \cite{cm}.  Here we just mention that,
in our \sy, a local stable manifold (LSM) $W^s(x)$ is a smooth curve
containing $x$ and whose main property is that, for all $y \in
W^s(x)$, $\lim_{n \to +\infty} \mathrm{dist} (T^n x,\, T^n y) = 0$,
where $\mathrm{dist}$ is the natural Riemannian distance in $\mM$
(with the convention that, if $x$ and $y$ belong to different
connected components of $\mM$, $\mathrm{dist}(x,y) = \infty$). A local
unstable manifold (LUM) $W^u(x)$ has the analogous property for the
limit $n \to -\infty$.

The \sy\ has a \hyp\ structure \emph{\`a la} Pesin, in the following
sense:

\begin{theorem}
  \label{thm-lsums}
  For $\mu$-a.e.\ $x \in \mM$ there is a LSM $W^s(x)$ and a LUM
  $W^u(x)$. The corresponding two foliations --- more correctly,
  laminations --- can be chosen invariant, namely $T W^s(x) \subset
  W^s(Tx)$ and $T^{-1} W^u(x) \subset W^u(T^{-1} x)$. Also, when
  endowed with a Lebesgue-equivalent $1$-dimensional transversal
  \me, they are absolutely continuous w.r.t.\ $\mu$.
\end{theorem}

The next theorem is the core technical result for all the proofs that
follow. It is not by chance that, in the field of \hyp\ \bi s, this is
called the \emph{fundamental theorem}.

\begin{theorem}
  \label{thm-connect}
  Given $n \in \Z$, $j \in \{1,2\}$ and a full-\me\ $A \subset
  \mN_n^j$, there exists a full-\me\ $B \subset \mN_n^j$ such that all
  pairs $x,y \in B$ are connected via a polyline of alternating LSUMs
  whose vertices lie in $A$. This means that, for $x,y \in B$, there
  is a finite collection of LSUMs, $W^s(x_1)$, $W^u(x_2)$, $W^s(x_3)$,
  $\dots$, $W^u(x_m)$, with $x_1 = x$, $x_m = y$, and such that each
  LSUM intersect the next transversally in a point of $A$.
\end{theorem}

The above theorems are proved in \cite{l2} for Lorentz gases that are
effectively two-dimensional and whose scatterers are smooth, i.e.,
$K=1$ in (A2). The first of the two differences is absolutely
inconsequential. The second affects the \emph{singularity set} of $T$,
that is, the set of all $x \in \mM$ whose \tr y, up to the next
crossing of a transparent wall, is singular. It is a well-known and
easily derivable fact that, in each component $\mN_n^j$ of the cross
section, the singularity set is a union of smooth curves, each of
which is associated to a specific source of singularity within the
cell $C_n$ (a tangential scattering, a vertex, the endpoint of a gate)
and an \emph{itinerary} of visited scatterers before that. Since both
the number of scatterers in each cell and the number of vertices per
scatterer are bounded, there can only be a finite number of
\emph{singularity lines} in each $\mN_n^j$. With this provision, the
proofs of \cite{l2} work in this case as well.

(In truth, the actual proofs are found in \cite{l1}, where the
existence of a \hyp\ structure and the fundamental theorem are shown
for the standard \bi\ cross section. In \cite{l2} these are extended
to the transparent cross section.  The idea behind the results of
\cite{l1} is this: Assumptions (A2)-(A4) guarantee that the geometric
features of the LT are ``uniformly good''. Then a refinement of a
standard trick ensures that most \o s of the \sy\ do not approach the
singularity set too fast, so that, in the construction of the \hyp\
structure, one can practically neglect them. As for the fundamental
theorem, all the local arguments in the classical proofs of Sinai and
followers for compact \bi s apply --- notice that we have uniform \hyp
ity and no \emph{cusps}, namely, zero-angle corners. The global
arguments have to do essentially with controlling the neighborhoods of
certain portions of the singularity set, which can be done with the
above-mentioned trick.)

\section{Recurrence}
\label{sec-rec}

We are interested in the recurrence and ergodic properties of the LTs
defined earlier.  To this goal, let us recall some definitions that
may not be obvious for infinite-\me\ \dsy s.

\begin{definition}
  \label{def-po-rec}
  The \me-preserving \dsy\ $(\mM, T, \mu)$ is called
  \emph{(Poincar\'e) recurrent} if, for every measurable $A \subseteq
  \mM$, the \o\ of $\mu$-a.e.\ $x \in A$ returns to $A$ at least once
  (and thus infinitely many times, due to the invariance of $\mu$).
\end{definition}

\begin{definition}
  \label{def-erg}
  The \me-preserving \dsy\ $(\mM, T, \mu)$ is called \emph{ergodic} if
  every $A \subseteq \mM$ measurable and invariant $\bmod \mu$ (that
  is, $\mu (T^{-1} A \, \symmdiff A) = 0)$, has either zero measure or
  full measure (that is, $\mu (\mM \setminus A) = 0$).
\end{definition}

If the \sy\ in question is an LT as introduced in Section
\ref{sec-pre} ($T = T_\l$ for some $\l \in \Omega^\Z$), it is proved
in \cite{l1, l2} that 
\begin{theorem}
 \label{thm-rec-iff-erg}
 $(\mM, T_\l, \mu)$ is \erg\ if and only if it is recurrent.
\end{theorem}

Understandably, proving recurrence (and thus ergodicity) of
\emph{every} \sy\ in the quenched random ensemble might be a daunting
task. It is possible, however, to prove it for a \emph{typical}
\sy. This will be achieved via a general result by Schmidt \cite{s} on
the recurrence of commutative cocycles over finite-\me\ \dsy s. We
state it momentarily.
\begin{definition}
  \label{def-co}
  Let $(\Sigma, F, \lambda)$ be a \pr-preserving dynamical \sy, and
  $f$ a measurable \fn\ $\Sigma \longrightarrow \Z^{\nu}$. The family
  of \fn s $\{ S_n \}_{n\in \N}$, defined by $S_0(\xi) \equiv 0$ and,
  for $n \ge 1$,
  \begin{displaymath}
    S_n(\xi) := \sum_{k=0}^{n-1} (f \circ F^k) (\xi)
  \end{displaymath}
  is called the \emph{cocycle of $f$}. Any such family is generically
  called \emph{commutative, $\nu$-dimensional, discrete cocycle}.
\end{definition}

\begin{theorem}
  \label{thm-co-rec}
  Assume that $(\Sigma, F, \lambda)$ is \erg\ and denote by $Q_n$ the
  distribution of $S_n / n^{1/\nu}$ relative to $\lambda$, i.e., the
  distribution on $\R^\nu$ defined by
  \begin{displaymath}
    Q_n (A) := \lambda \left( \rset{\xi \in \Sigma} {\frac{S_n 
    (\xi)} {n^{1/\nu}} \in A} \right),
  \end{displaymath}
  where $A$ is any measurable set of $\R^\nu$. If there exists a
  positive-density sequence $\{ n_k \}_{k\in\N}$ and a constant
  $\kappa > 0$ such that
  \begin{displaymath}
    Q_{n_k}( \mathcal{B}(0,\rho) ) \ge \kappa \rho^{\nu} 
  \end{displaymath}
  for all sufficiently small balls $\mathcal{B}(0,\rho) \subset
  \R^\nu$ (of center 0 and radius $\rho$), then the cocycle $\{ S_n
  \}$ is \emph{recurrent}, namely, for $\lambda$-a.e.\ $\xi \in
  \Sigma$, there exists a subsequence $\{ n_j \}_{j \in \N}$ 
  such that
  \begin{displaymath}
    S_{n_j}(\xi) = 0, \quad \forall j \in \N.
  \end{displaymath}  
\end{theorem}

The above result is a slight weakening of the original theorem by
Schmidt, whose proof can be found in \cite{s}. (In truth, the original
formulation required $F$ to be invertible mod $\lambda$. The
generalization to non-invertible \me-preserving maps is an easy
exercise which can be found, e.g., in \cite[App.~A.2]{l3}).

\skippar

In the following we will introduce a suitable \pr-preserving \dsy\ and
a 1-dimensional cocycle with the property that the recurrence of the
latter is equivalent to the Poincar\'e recurrence of $\Pi$-a.e.\ LT
$\l$ (we call this situation \emph{almost sure recurrence} of the
quenched random LT; details in Section \ref{sec-povp}). Observe that,
for $\nu = 1$, the quantity $S_n / n^{1/\nu}$ is precisely the
Birkhoff average of $f$. Thus the \erg ity of $(\Sigma, F, \lambda)$,
which implies the law of large numbers for $\{ S_n \}$, is enough to
apply Theorem \ref{thm-co-rec}.

\section{The point of view of the particle}
\label{sec-povp}

For $j \in \{ 1,2 \}$, let us consider $\mN_0^j$ as defined in
(\ref{def-nnj}), and rename it $\mN^j$ for short.  In this section we
will work extensively with the cross-section $\mN := \mN^1 \cup
\mN^2$.

Let us call $\mu_0$ the standard \bi\ \me\ on $\mN$, normalized to
$1$. If $\om\in\Omega$ determines the configuration of scatterers in
$C_0$, we can define a map $R_\om: \mN \longrightarrow \mN$ as follows
(cf.\ Fig.~\ref{imgromega}).  Trace the forward \tr y of $x := (q,v)
\in \mN$ until it crosses $G^1$ or $G^2$ for the first time (almost
all \tr ies do). This occurs at a point $q_1$ with velocity $v_1$. If,
for $\epsilon \in \{-1,+1\}$, $C_{\epsilon}$ is the cell that the
particle enters upon leaving $C_0$, define
\begin{eqnarray}
  \label{eq-def-R}
  R_\om \, x = R_\om (q,v) &:=& (\tau^{-\epsilon} (q_1), v_1) \in \mN, \\
  \label{eq-def-e}
  e(x,\om) &:=& \epsilon.
\end{eqnarray}
We name $e$ the \emph{exit function}. From our earlier discussion on
the transparent cross sections, $R_\om$ preserves $\mu_0$.

\fig{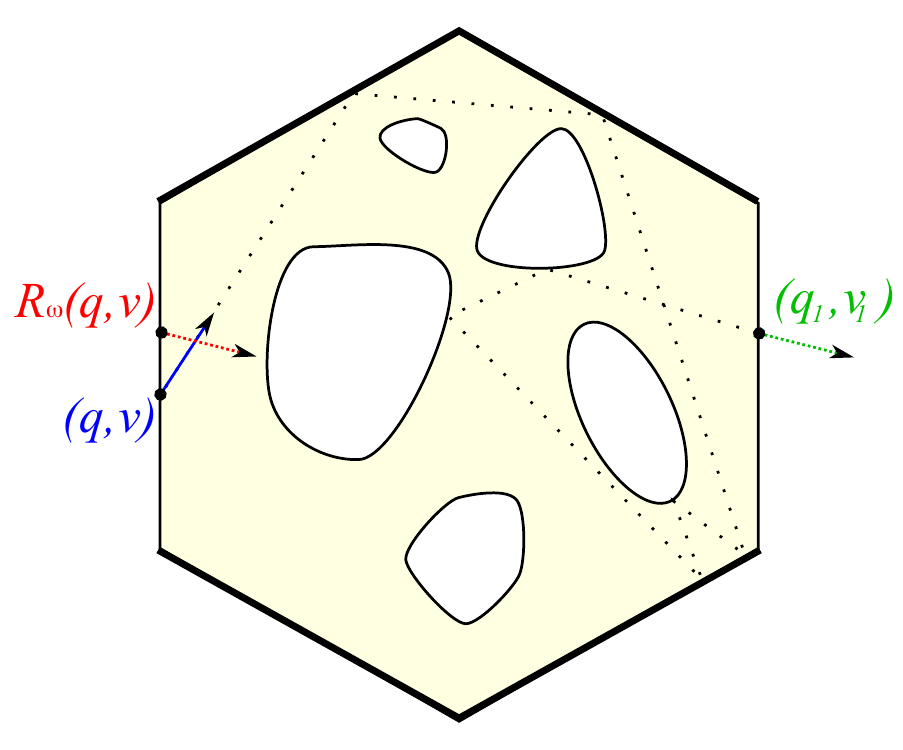}{0.5\linewidth}{The definition of the map $R_{\omega}$.}

We introduce the \dsy\ $(\Sigma, F, \lambda)$, where
\begin{itemize}
\item $\Sigma := \mN\times\Omega^{\Z}$.

\item $F(x,\l) := (R_{\l_0} x, \sigma^{e(x,\l_0)} (\l))$, defining a
  map $\Sigma \longrightarrow \Sigma$. Here $\l_0$ is the $0$th
  component of $\l$ and $\sigma$ is the left shift on $\Omega^{\Z}$,
  introduced in (A1) (therefore $\sigma^\epsilon (\l) =
  \{\l'_n\}_{n\in\Z}$, with $\l'_n := \l_{n + \epsilon}$).

\item $\lambda := \mu_0 \times \Pi$. Clearly, $\lambda (\Sigma) =
  1$. Also, using that $F$ is invertible, $R_\om$ preserves $\mu_0$
  for every $\om \in \Omega$, and $\sigma$ preserves $\Pi$, it can be
  seen that $F$ preserves $\lambda$. (This is ultimately a consequence
  of the fact that every LT preserves the same \me.)
\end{itemize}

The idea behind this definition is that, instead of following a given
\o\ from one cell to another, we every time shift the LT in the
direction opposite to the \o\ displacement, so that the point always
lands in $C_0$. For this reason the \dsy\ just introduced is called
\emph{the point of view of the particle}.  Clearly, $F: \Sigma
\longrightarrow \Sigma$ encompasses the dynamics of all points on all
realizations of $\Omega^{\Z}$.

\begin{proposition}
  \label{prop-rec-co}
  If the cocycle of the exit \fn\ $e$ is recurrent, then the quenched
  random LT is \emph{almost surely recurrent} in the sense that, for
  $\Pi$-a.e.\ $\l \in \Omega^\Z$, $(\mM, T_\l, \mu)$ is recurrent.
\end{proposition}

\proof Before starting the actual proof, we recall that an easy
argument \cite[Prop.~2.6]{l2} shows that the extended \sy\ $(\mM,
T_\l, \mu)$ is either recurrent or totally dissipative (i.e.,
transient): no mixed situations occur. Therefore, the existence of one
\emph{recurring set} (i.e., a positive-\me\ set $A$ such that
$\mu$-a.a.\ points of $A$ return there at some time in the future) is
enough to establish the same property for \emph{all} measurable sets.

Now, calling $\{ S_n \}$ the cocycle of $e$, the hypothesis of
Proposition \ref{prop-rec-co} amounts to saying that, for
$\lambda$-a.e.\ $(x,\l) \in \Sigma$, there exists $n = n(x,\l)$ such
that $S_n(x,\l) = 0$. That is, considering the LT $\l$, $T_\l^n x \in
\mN_0$ (recall that $x \in \mN_0$ by construction). Let us call such a
pair $(x,\l)$ \emph{typical}.

By Fubini's Theorem, $\Pi$-a.a.\ $\l \in \Omega^\Z$ are such that
$(x,\l)$ is typical for $\mu_0$-a.a.\ $x \in \mN$. For such $\l$,
$\mN_0 = \mN$ is a recurring set of $T_\l$, therefore $(\mM, T_\l,
\mu)$ is recurrent.
\qed

As it was mentioned at the end of Section \ref{sec-rec}, the
recurrence of the cocycle of $e$ is implied by ergodicity of $(\Sigma,
F, \lambda)$. On the other hand,

\begin{theorem}
  \label{thm-povp-erg}
  Under assumptions \emph{(A1)-(A5)}, the \dsy\ $(\Sigma, F, \lambda)$
  defined above is \erg.
\end{theorem}

\proof The proof can be divided in three steps:

\begin{enumerate}
\item Every ergodic component of $(\Sigma, F,\lambda)$ is of the form
  $\bigcup_{j=1}^2 \mN^{j} \times B_j$ $\bmod \lambda$, where $B_j$ is
  a measurable set of $\Omega^\Z$.

\item $\Pi(B_j) \in \{ 0,1 \}$.

\item There is only one \erg\ component.
\end{enumerate}
We now describe each step separately.
\begin{enumerate} 
\item For a fixed $\l$, consider the extended \dsy\ $(\mM, T_\l,
  \mu)$, for which Theorem \ref{thm-lsums} holds. Through the obvious
  isomorphism, copy those LSUMs of the extended \sy\ which are
  included in $\mN_0$ onto $\mN \times \{\l\}$. These may be called
  \emph{LSUMs for the fiber} $\mN \times \{\l\}$ (although $(\Sigma,
  F, \lambda)$ cannot be regarded as a \emph{bona fide} \hyp\
  \dsy). By Theorem \ref{thm-connect}, in each connected component of
  $\mN \times \{\l\}$, namely, $\mN^1 \times \{\l\}$ and $\mN^2 \times
  \{\l\}$, a.e.\ pair of points can be connected through a sequence of
  LSUMs for the fiber, intersecting at typical points. Hence, via the
  usual Hopf argument \cite{cm}, the whole $\mN^j \times \{\l\}$ lies
  the same \erg\ component, at least for a.e.\ $\l$. Therefore an
  $F$-invariant set in $\Sigma$ can only come in the form $I =
  \bigcup_{j=1}^2 \mN^j \times B_j$. That $B_j$ is measurable is a
  consequence of Lemma A.1 in \cite{l2}.

\item If $I$ as written above is $F$-invariant, then $\mN^1 \times
  B_1$ is $F_1$-invariant, where $F_1$ is the first-return map of $F$
  onto $\mN^1 \times \Omega^\Z$. Consider a typical $\l \in B_1$ in
  the following sense: for $\mu_0$-a.e.\ $x \in \mN^1$, the
  $F_1$-orbit of $(x,\l)$ is entirely included in $\mN^1 \times B_1$;
  also, looking at (A5), the LT $\l$ possesses a positive-\me\ set of
  \tr ies entering $C_0$ through $G^1$ and leaving it through
  $G^2$. This implies that there exists an $x \in \mN^1$ such that
  $F(x,\l) \in \mN^1 \times B_1$ and $F(x,\l) = (x',\sigma(\l))$, for
  some $x'$.  Hence $\sigma(\l) \in B_1$. Considering that this
  happens for $\Pi$-a.a.\ $\l \in B_1$, we obtain $\sigma(B_1)
  \subseteq B_1$ $\bmod\ \Pi$. (A1) then implies that $\Pi(B_1) \in \{
  0,1 \}$. The analogous assertion for $B_2$ can be proved by using
  $F_2$, the first-return map onto $\mN^2 \times \Omega^\Z$; the
  existence of a non-singular \tr y going from $G^2$ to $G^1$, and
  $\sigma^{-1}$ instead of $\sigma$.

\item It cannot happen that $\mN^1 \times \Omega^\Z$ and $\mN^2 \times
  \Omega^\Z$ are two different \erg\ components, because, via (A5),
  for $\Pi$-a.e.\ $\l \in \Omega^\Z$ there is a positive $\mu_0$-\me\
  of points $x \in \mN^1$ for which $F(x,\l) \in \mN^2 \times
  \Omega^\Z$.
\end{enumerate}
\vspace{-6pt}
\qed

As explained in the last paragraph of Section \ref{sec-rec},
Proposition \ref{prop-rec-co} and Theorem \ref{thm-povp-erg} yield our
main result:

\begin{theorem}
  \label{thm-main}
  Under assumptions \emph{(A1)-(A5)}, $(\mM, T_\l, \mu)$ is recurrent
  for $\Pi$-a.e.\ $\l \in \Omega^\Z$.
\end{theorem}

\section{Extensions}
\label{sec-ext}

If we look at the proof of Theorem \ref{thm-povp-erg}, it is apparent
that its key argument is that each horizontal fiber $\mN^j \times
\Omega^\Z$ is part of the same \erg\ component. Once that is known,
one simply uses (A5) to show that a given ergodic component invades
the whole phase space, first for the map $F_j$ and then for the map
$F$ itself. The details of the dynamics are not relevant for this
argument.

By Theorem \ref{thm-co-rec}, the ergodicity of the point of view of
the particle implies the recurrence of our cocycle, because the
cocycle is one-dimensional. Thus, as long as we deal with systems in
which the position of the particle can be described, in a discrete
sense, by a one-dimensional cocycle, the foregoing arguments can be
used to prove the almost sure recurrence of a more general class of
LTs.

In the present section we sketch the construction of some of these
extensions.

\subsubsection*{Same gates, different cells}

There is no reason why all the cells $C_n$ should be the same
polygon. One can easily consider random cells $C_n$ in which the
border too depends on the random parameter $\l_n$. This can be devised
by putting extra flat scatterers in a sufficiently large cell in order
to produce any desired shape; see Fig.~\ref{imgcella}.  As long as
each cell has two opposite congruent gates and (A1)-(A5) are verified,
all the previous results continue to hold.

\fig{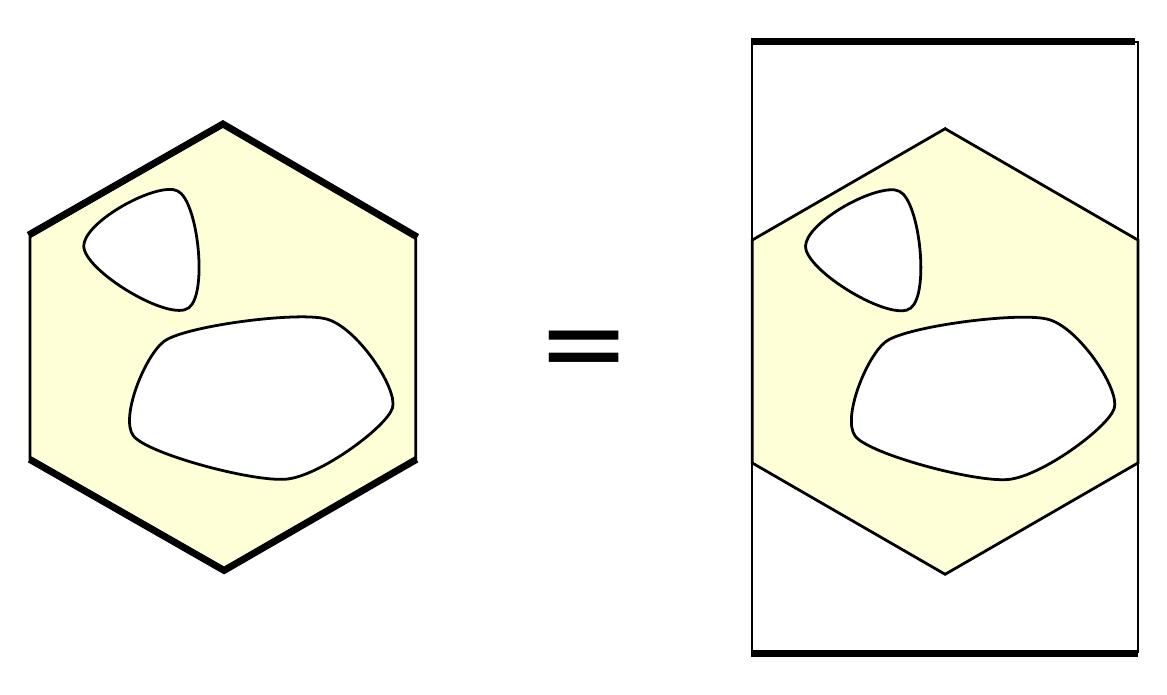}{0.5\linewidth}{Realizing a randomly-shaped cell out of
  a standard cell.}

In fact, one can allow for the distance between the gates to vary with
$\l_n$ as well (in (\ref{eq-def-R}) simply replace $\tau^{-\epsilon}$
with the cell-dependent local translation
$\tau^{-\epsilon}_{\omega}$). An example of this type of LT is shown
in Fig.~\ref{imgtubodiverso}.

\fig{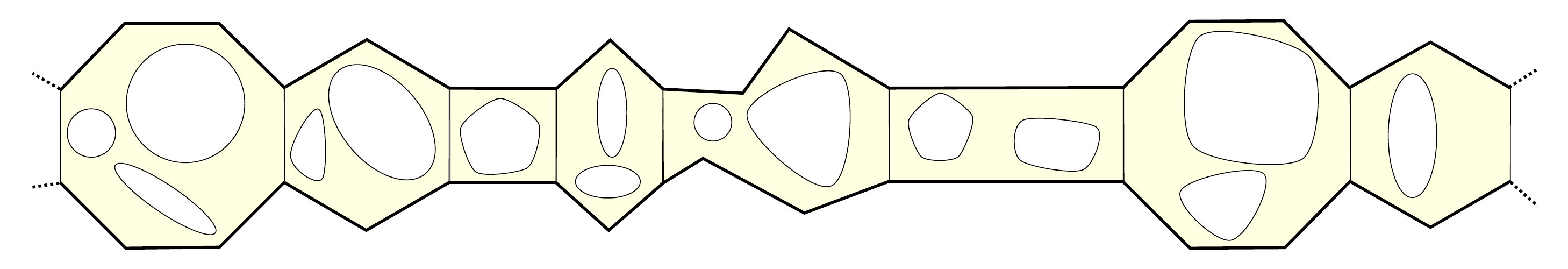}{0.9\linewidth}{An LT with different cells.}

\subsubsection*{Same cells, poly-gates}

One can also define $G^j$ to be the union of a finite number of sides
$G^{ji}$, with $i$ varying in some index set $I$, provided that there
is a translation $\tau$ such that $\tau(G^1) = G^2$; see
Fig.~\ref{imgesagoni2}. However, in order for steps 2 and 3 of the
proof of Theorem \ref{thm-povp-erg} to hold, (A5) needs to be replaced
by
\begin{itemize}
\item[(A5')] For a.e.\ $\ell$, all $j,j' \in \{ 1,2 \}$ and all $i,i'
  \in I$, there is a non-singular \tr y entering $C_0$ through
  $G^{ji}$ and leaving it through $G^{j'i'}$.
\end{itemize}
\fig{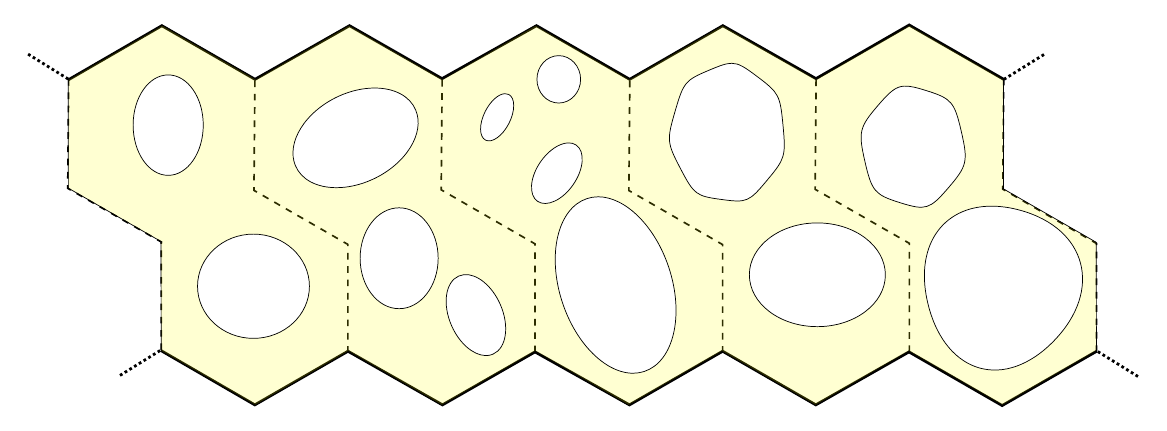}{0.7\linewidth}{An LT with non-trivial gates.}

\subsubsection*{From translation to general isometry}

Another hypothesis that is not crucial is that $G^1$ is mapped onto
$G^2$ via a translation. One can imagine that $\Z$ acts upon the
Lorentz tube via a general isometry, for example a roto-translation,
as in Fig.~\ref{imgtuboarrotolato1}.

\fig{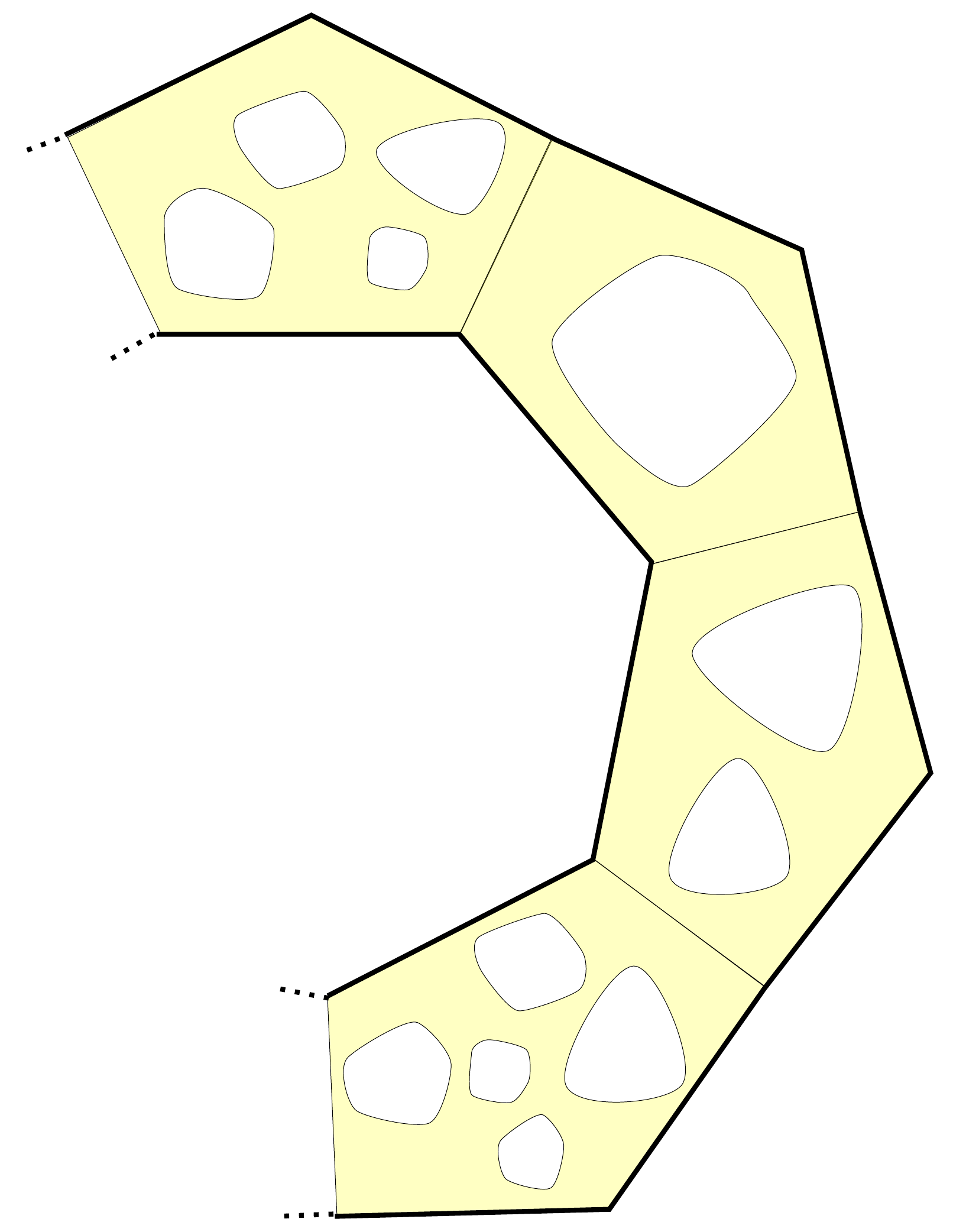}{0.35\linewidth}{A spiraling LT.}

The only problem, in this case, is that, quite generally, the
resulting tube will have self-intersections. One can simply do away
with it by disregarding the self-intersections, e.g., by declaring
that any two portions of the tube that intersect in the plane actually
belong to different sheets of a Riemann surface.

\subsubsection*{Random gates and random isometries}

Assume that the fundamental domain is a polygon $C$ such that $p$ of
its sides ($p \ge 2$) are congruent. In this case it is possible to
randomize the choice of the gates too. That is, one can let the random
parameter $\l_n$ decide which of the $p$ congruent sides of $C_n$
will play the role of the ``left'' and ``right'' gates.
Moreover, $\l_n$ can also prescribe how the right gate of $C_n$
attaches to the left gate of $C_{n+1}$; see Fig.~\ref{imgtubonuovo}.

\fig{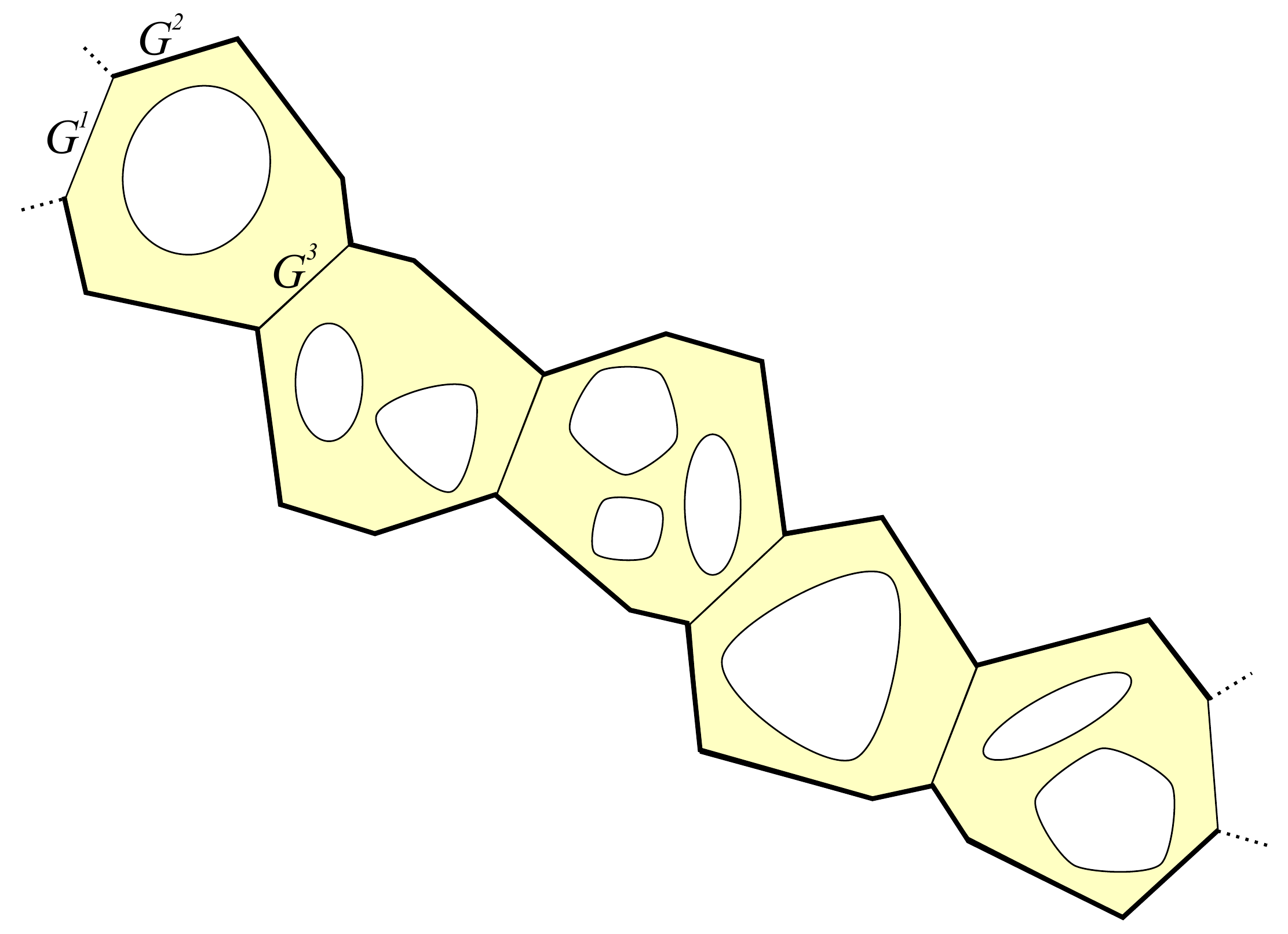}{0.65\linewidth}{An LT with random gates (in this
  case $p=3$, see text).}

In order to implement this idea, we need to slightly change our
previous notation. Let $\{ G^j \}_{j=1}^p$ be a fixed ordering of the
$p$ congruent sides of $C$ mentioned above. For any such $j$, let
$\mN^j$ denote the transparent, incoming, cross section relative to
$G^j$, as in (\ref{def-nnj}). Then set $\mN := \bigcup_j \mN^j$.

We assume that there exist two \fn s $j_1, j_2: \Omega \longrightarrow
\{ 1, \ldots, p \}$ such that $j_1(\om) \ne j_2(\om)$, $\forall
\om$. This is how $\om$ specifies that $G^{j_1}$ and $G^{j_2}$ are the
left and right gates, respectively, of $C$.

In lieu of $R_\om$, cf.\ (\ref{eq-def-R}), we use the more general map
$R_\l: \mN \longrightarrow \mN$ defined as follows. For $x = (q,v) \in
\mN$, let $G^j$ be the first side of its kind that the forward
flow-\tr y of $x$ hits within $C$, and denote by $q_1$ and $v_1$,
respectively, the hitting point in $G^j$ and the precollisional
velocity there (see Fig.~\ref{imgromega}).
\begin{itemize}

\item If $j = j_2(\l_0)$ then $R_\l \, x := \xi_{\l_0} \circ
  \rho_{j_2(\l_0), j_1(\l_1)} (q_1, v_1)$. Here $\rho_{j, j'}$ is the
  transformation that rigidly maps the outer pairs $(q_1,v_1)$ based
  in $G^j$ onto the inner pairs based in $G^{j'}$ (it is a
  rototranslation in the $q$ variable); and $\xi_\om: \mN
  \longrightarrow \mN$, depending on the usual random parameter $\om$,
  is either the identity or the transformation that flips all the
  segments $G_j$ and changes the $v$ variable accordingly. So, through
  $\xi_\om$, $\l_n$ decides whether $C_n$ and $C_{n+1}$ have the same
  or opposite orientations (cf.\ Fig.~\ref{imgtubonuovo}).  In this
  case, the exit \fn\ is set to the value $e(x, \l_0) := 1$.

\item If $j = j_1(\l_0)$ then, in accordance with the previous case,
  $R_\l \, x := \xi_{\l_{-1}} \circ \rho_{j_1(\l_0), j_2(\l_{-1})}
  (q_1, v_1)$ (notice that $\xi_\om^{-1} = \xi_\om$). In this case,
  $e(x, \l_0) := -1$.

\item For all the other $j$, $R_\l \, x := (q_1, v_2)$, where $v_2 :=
  v_1 + 2 (v_1 \cdot o_j) o_j$ is the postcollisional velocity
  corresponding to a \bi\ bounce against $G^j$ with incoming velocity
  $v_1$ ($o_j$ denoted the inner normal to $G^j$). For this last case,
  $e(x, \l_0) := 0$.
\end{itemize}

\end{document}